\documentclass[12pt]{article}
\usepackage{amssymb}
\usepackage{amsmath,amsthm}
\usepackage[latin1]{inputenc}
\usepackage{graphicx}
\usepackage{enumerate}
\usepackage{tkz-graph}
\usepackage{pstricks}
\usepackage{tikz}
\usepackage{verbatim}

\newtheorem{theorem}{Theorem}

\newtheorem{corollary}[theorem]{Corollary}

\newtheorem{lemma}[theorem]{Lemma}

\newtheorem{proposition}[theorem]{Proposition}
\newtheorem{remark}[theorem]{Remark}

\title{Computing the local metric dimension of a graph from the local metric dimension of primary  subgraphs}

\author{Juan A. Rodr\'{i}guez-Vel\'{a}zquez, Carlos Garc\'{\i}a G\'omez and \\
Gabriel A. Barrag\'{a}n-Ram\'{\i}rez
    \\
{\small Departament d'Enginyeria Inform\`atica i Matem\`atiques,}\\
{\small Universitat Rovira i Virgili,}  {\small Av. Pa\"{\i}sos
Catalans 26, 43007 Tarragona, Spain.} \\{\small
juanalberto.rodriguez\@@urv.cat, gbrbcn\@@gmail.com}
}

\begin{document}
\maketitle

\begin{abstract}
For an ordered subset $W = \{w_1, w_2,\dots w_k\}$ of vertices and a vertex $u$ in a connected graph $G$, the representation of $u$ with respect to $W$ is the ordered $k$-tuple $ r(u|W)=(d(v,w_1), d(v,w_2),\dots,$ $d(v,w_k))$, where $d(x,y)$ represents the distance between the vertices $x$ and $y$. The set $W$ is a local metric generator for $G$ if every two adjacent vertices of $G$ have distinct representations. A minimum local metric generator is called a \emph{local metric
basis} for $G$ and its cardinality the \emph{local metric dimension} of G.   We show that the computation of the local metric dimension of a graph with cut vertices is reduced to the computation  of the local metric dimension of the so-called primary subgraphs.
The main results are applied to specific constructions including bouquets of graphs,  rooted product graphs, corona product graphs, block graphs and chain of graphs.
\end{abstract}

\section{Introduction}

 A {\em generator} of a metric space is a set $S$ of points in the space with the property that every point of the space is uniquely determined by its distances from the elements of $S$. Given a simple and connected graph $G=(V,E)$, we consider the metric $d_G:V\times V\rightarrow \mathbb{N}$, where $d_G(x,y)$ is the length of a shortest path between $x$ and $y$. $(V,d_G)$ is clearly a metric space. A \emph{metric generator} of a
connected graph $G$ is a subset of vertices, $W\subset V(G),$ for which,
given any pair of vertices $u,v\in V(G)$ there is at least one element  $w\in W$
 for which we have%
\[
d_{G}(u,w)\neq d_{G}(v,w).
\]
We say then, that $w$ is able to distinguish   the pair of vertices $u,v.$ A metric generator
 with minimum cardinality  is defined
as a \emph{metric basis} for $G.$ The cardinality of this set is denoted by
$\dim(G)$ and is referred as the metric dimension of $G.$

We can see a metric basis $S$ of $G$ as an ordered set  $S=\{s_{1}, s_{2}, \ldots, s_{d}\}$. In this sense,  we refer to the vector
$$r(u|S)=\left(d_{G}(u,s_{1}),d_{G}(u,s_{2}),\ldots,d_{G}(u,s_{d})\right)$$
as the  coordinate vector of $u$ with respect to the basis $S$. Note that since $S$ is a   metric basis  for $G$, for any pair of vertices $u$ and $v$ of $G$, it holds that $r(u|S)\neq r(v|S)$. Hence, each vertex is uniquely determined by its  coordinate vector with respect to a basis.

These concepts were first
introduced by Slater in \cite{Slater1975}, where the metric generators were called \emph{locating sets}. The concept of metric dimension of a connected graph was
introduced independently by Harary and Metler in \cite{Harary1976}, where metric generators received the name of \emph{resolving sets.}
 After these papers were published several authors developed diverse theoretical works about this topic, for instance, we cite \cite{Caceres2007, Hernando2005, Chappell2008, Chartrand2000, Chartrand2001a, Chartrand2003b, Chartrand2000a, Chartrand2000b, Enomoto2000, Enomoto2002, Fehr2006, Harary1976, Haynes2006, Khuller1996, Khuller1994, Nakamigawa2003, Saenpholphat2003a, Saenpholphat2003b, Saenpholphat2003c}.
 Slater described the usefulness of these ideas into long range aids to navigation \cite{Slater1975}. Also, these concepts  have some applications in chemistry for representing chemical compounds \cite{Johnson1998,Johnson1993} or to problems of pattern recognition and image processing, some of which involve the use of hierarchical data structures \cite{Melter1984}. Other applications of this concept to navigation of robots in networks and other areas appear in \cite{Chartrand2000,Hulme1984,Khuller1996}.

In this paper we are interested in a local version of  metric generators introduced by Okamoto et al. in \cite{Okamoto2010}.  Given a
connected graph $G,$ we   define a \emph{local metric generator} as a set of
vertices that  distinguishes  any pair of  adjacent vertices. This means that, given
two adjacent vertices $u,v\in V(G)$  there is at least an element of this set, say $w$, for which we have $
d_{G}(u,w)\neq d_{G}(v,w)$.



If a local metric generator has minimum  cardinality among all  local metric generators, then we call this set
a \emph{local metric basis} for $G.$ The cardinality of the local metric basis is
denoted by $\dim _{l}(G)$ and it is called the \emph{local metric dimension} of $G$.
Note that each metric generator is also a local metric
generator because each metric generator distinguishes any pair of vertices,  while a local metric generator
only distinguishes pairs of neighbours. Then the following relation between the local metric
dimension and the metric dimension of a graph is valid
\[
1\leq \dim_{l}(G)\leq \dim (G)\leq n-1.
\]

In this paper we show that the computation of the local metric dimension of a graph with cut vertices is reduced to the computation  of the local metric dimension of the so-called primary subgraphs.
The main results are applied to specific constructions including bouquets of graphs,  rooted product graphs, corona product graphs, block graphs and chain of graphs.

The following basic results, established in 
 \cite{Okamoto2010},   will be used in this paper.

\begin{theorem}\label{ldimb01}
{\rm \cite{Okamoto2010}} Let $G$ be a nontrivial connected graph of order $n.$ Then $\dim
_{l}(G)=n-1$ if and only if $G=K_{n}$ and $\dim _{l}(G)=1$ if and only if $G$
is bipartite.
\end{theorem}

\begin{theorem}\label{ldimb02}
{\rm \cite{Okamoto2010}}. A connected graph $G$ of order $n\geq 3$ has local metric dimension
$\dim _{l}(G)=n-2$ if and only if the clique number  of $G$ is $
\omega (G)=n-1.$
\end{theorem}

In this work the remain definitions are given the first time that the concept is found in
the text.

\section{Main results}

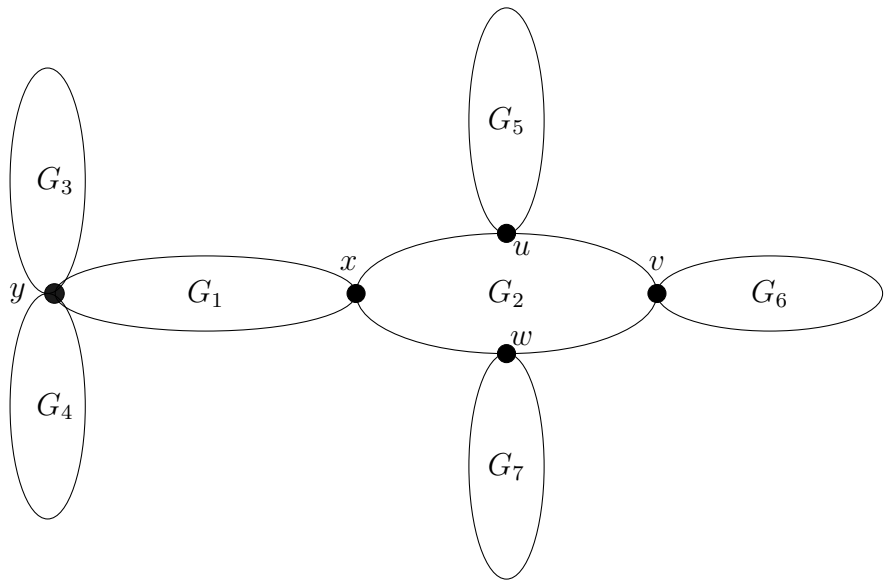
\begin{figure}[h]
   \centering
\begin{tikzpicture}
\draw  (-2,0) ellipse (2cm and 0.5cm);
\filldraw[fill opacity=1,fill=black]  (0,0) circle (0.12cm);
\draw  (2,0) ellipse (2cm and 0.8cm);
\draw  (2,2.3) ellipse (0.5cm and 1.5cm);
\filldraw[fill opacity=1,fill=black]  (2,0.8) circle (0.12cm);
\draw  (5.5,0) ellipse (1.5cm and 0.5cm);
\filldraw[fill opacity=1,fill=black]  (4,0) circle (0.12cm);
\draw (-4.1,1.5) ellipse (0.5cm and 1.5cm);
\draw (-4.1,-1.5) ellipse (0.5cm and 1.5cm);
\filldraw[fill opacity=0.9,fill=black]  (-4.01,0) circle (0.13cm);
\draw  (2,-2.3) ellipse (0.5cm and 1.5cm);
\filldraw[fill opacity=1,fill=black]  (2,-0.8) circle (0.12cm);
\node at (-2,0) {$G_1$ }; 
\node at (-4.5,0) {$y$ };
\node at (2,0) {$G_2$ };
\node at (-0.1,0.4) {$x$ };
\node at (2,2.3) {$G_5$ };
\node at (2,-2.3) {$G_7$ };
\node at (-4,1.5) {$G_3$ };
\node at (-4,-1.5) {$G_4$ };
\node at (5.5,0) {$G_6$ };
\node at (4,0.4) {$v$ };
\node at (2.2,0.6) {$u$ };
\node at (2.2,-0.6) {$w$ };
\end{tikzpicture}
\caption{A graph $G[{\cal H}]$ obtained by point-attaching from ${\cal H}=\{G_1,G_2, . . . ,G_7\}$}
\label{point-attaching}
\end{figure}

Let $G[{\cal H}]$ be a connected graph constructed from a family of pairwise disjoint (non-trivial) connected graphs
${\cal H}=\{G_1, . . . ,G_k\}$ as follows. Select a vertex of $G_1$, a vertex of $G_2$, and identify these two
vertices. Then continue in this manner inductively. More precisely, suppose that we
have already used $G_1, . . . ,G_i$ in the construction, where $2\le i \le k- 1$. Then select
a vertex in the already constructed graph (which may in particular be one of the
already selected vertices) and a vertex of $G_{i+1}$; we identify these two vertices. Note
that any graph $G[{\cal H}]$ constructed in this way has a tree-like structure, the $G_i's$ being its
building stones (see Figure \ref{point-attaching}). 
  
  We will briefly say that $G[{\cal H}]$ is obtained by \textit{point-attaching}
from $G_1, . . . ,G_k$ and  that $G_i's$ are the \textit{primary subgraphs} of $G[{\cal H}]$. We will also say that the vertices of $G[{\cal H}]$ obtained by identifying two vertices of different primary subgraphs are the \textit{attachment vertices} of $G[{\cal H}]$.  The above terminology was previously introduced in  \cite{Deutsch-Klavzar2013} where the authors obtained an expression that reduces the computation of the
Hosoya polynomials of a graph with cut vertices to the Hosoya polynomial of
the so-called primary subgraphs.

     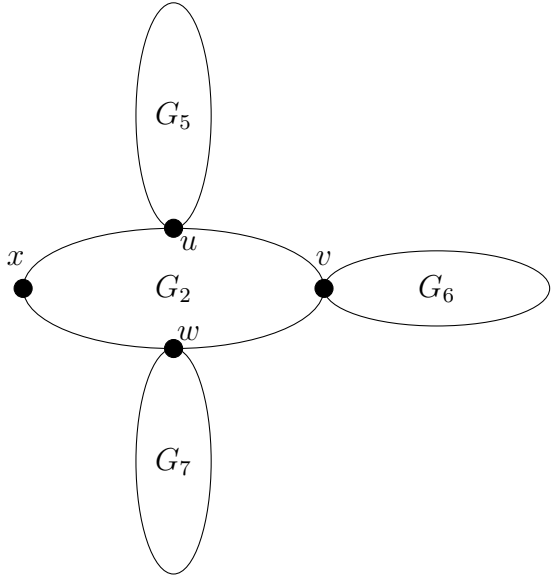
\begin{figure}[h]
   \centering
\begin{tikzpicture}
\filldraw[fill opacity=1,fill=black]  (0,0) circle (0.12cm);
\draw  (2,0) ellipse (2cm and 0.8cm);
\draw  (2,2.3) ellipse (0.5cm and 1.5cm);
\filldraw[fill opacity=1,fill=black]  (2,0.8) circle (0.12cm);
\draw  (5.5,0) ellipse (1.5cm and 0.5cm);
\filldraw[fill opacity=1,fill=black]  (4,0) circle (0.12cm);
\draw  (2,-2.3) ellipse (0.5cm and 1.5cm);
\filldraw[fill opacity=1,fill=black]  (2,-0.8) circle (0.12cm);

\node at (2,0) {$G_2$ };
\node at (-0.1,0.4) {$x$ };
\node at (2,2.3) {$G_5$ };
\node at (2,-2.3) {$G_7$ };

\node at (5.5,0) {$G_6$ };
\node at (4,0.4) {$v$ };
\node at (2.2,0.6) {$u$ };
\node at (2.2,-0.6) {$w$ };
\end{tikzpicture}
\caption{The subgraph $G_1(x^+)$ of the graph $G[{\cal H}]$ shown in Figure \ref{point-attaching}.}
\label{Subgraphpoint-attaching}
\end{figure}
  
  To begin with the study of the local metric dimension of $G[{\cal H}]$ we need some additional terminology. Given an attachment  vertex $x$ of $G[{\cal H}]$ and a  primary subgraph $G_j$ such that $x\in V(G_j)$, we define the  subgraph $G_j(x^+)$  of $G[{\cal H}]$ as follows. We remove from $G[{\cal H}]$ all the edges connecting $x$ with vertices in  $G_j$, then $G_j(x^+)$  is the connected component which has $x$ as a vertex.
 For instance, Figure \ref{Subgraphpoint-attaching} shows the subgraph $G_1(x^+)$ of the graph $G[{\cal H}]$ shown in Figure \ref{point-attaching}. 
  
  Let $J_{\cal H}\subseteq [k]$ be the set of subscripts such that $j\in J_{\cal H}$ whenever 
  $G_j$ is a  non-bipartite primary subgraph of $G[{\cal H}]$. Note that $J_{\cal H}=\emptyset$  if and only if $G[{\cal H}]$ is bipartite, \textit{i.e.}, $J_{\cal H}=\emptyset$ if and only if $\dim_l(G[{\cal H}])=1$. From now on we assume that $  J_{\cal H} \ne \emptyset$.  

\newpage

Now, let 
 $C_j$ be the set composed by attachment  vertices of  $G[{\cal H}]$ belonging to $V(G_j)$ such that $x\in C_j$ whenever $G_j(x^+)$ is not bipartite. For instance, if $G_2$, $G_3$ and $G_7$ are the non-bipartite primary subgraphs of the graph shown in Figure \ref{point-attaching}, then $C_2=\{x,w\}$.

For any $j\in J_{\cal H}$ we define 
  $$\alpha_j=\max_{B\in {\cal B}(G_j)} \left\{ \vert C_j\cap B \vert \right\},$$
 where ${\cal B}(G_j)$ is the set of local metric bases of $G_j$, 
   \textit{i.e.}, $\alpha_j$ is the maximum cardinality of a set $ \{x_{j_1}, x_{j_2}, ...,x_{j_{\alpha_j}}\}\subseteq V(G_j)$  composed by  attachment  vertices of $G[{\cal H}]$ belonging simultaneously to a local metric basis of $G_j$ such that for every  $l\in \{1,...,\alpha_j\}$ the subgraph $G_j(x_{j_{l}}^+)$ is not  bipartite.

\begin{theorem}\label{SecondTheoremSeparableGraphs}
For any   non-bipartite graph $G[{\cal H}]$ obtained by point-attaching from a family of connected graphs  ${\cal H}=\{G_1, . . . ,G_k\}$,
$$\dim_l(G[{\cal H}])\le \sum_{j\in J_{\cal H}}(\dim_l(G_j)-\alpha_j).$$
\end{theorem}

\begin{proof}
For any $j\in J_{\cal H}$ we take  $B_j\in {\cal B}(G_j)$ and $M_j\subseteq B_j\cap C_j$ such that $\vert M_j\vert =\alpha_j$. We claim that 
$B=\displaystyle\bigcup_{j\in J_{\cal H}}(B_j-M_j)$ is a local metric generator for $G[{\cal H}]$. 

First of all, note that by the structure of $G[{\cal H}]$ we have that for any $v\in M_j$  there exists a non-bipartite primary subgraph
 $G_r$, which is a subgraph of  $G_j(v^+)$, such that $B_r-M_r\ne \emptyset$. To see this  we 
  take a non-bipartite primary subgraph  $G_{j_1}$, which is a subgraph of $G_j(v^+)$, next, if $B_{j_1}=M_{j_1}$, then we take $v_1\in V(G_{j_1})$ and, as above, we take a non-bipartite primary subgraph  $G_{j_2}$, which is a subgraph of $G_j(v_{1}^+)$, and if $B_{j_2}=M_{j_2}$ then we repeat this process until obtain  a non-bipartite primary subgraph  $G_{j_t}$, which is a subgraph of 
    $G_j(v_{t-1}^+)$ such that  $\vert B_{j_t}\vert >\vert M_{j_t}\vert$ (at worst, we will arrive to a subgraph $G_j(v_{t-1}^+)$ containing only one non-bipartite primary subgraph).
 With this fact in mind, we differentiate the following cases for two adjacent vertices $x,y\in V(G_i)$.
\\
\\
\noindent Case 1. $i\in J_{\cal H}$. If the pair $x,y$ is distinguished by some $u\in B_i-M_i$, then we are done. Now, if the pair $x,y$ is distinguished by $v\in  M_i$, 
 then we take $G_r$ as a non-bipartite  primary subgraph of  $G_i(v^+)$ such that $B_r-M_r\ne \emptyset$. Since
 the pair $x,y$ is distinguished by any vertex of $G_i(v^+)$, it is also distinguished by any $u\in B_r-M_r$.
 \\
\\
\noindent Case 2. $i\in [k]-J_{\cal H}$. In this case, we take $j\in J_{\cal H}$ such that $B_j-M_j\ne \emptyset$ and, since $G_i$ is bipartite, the pair $x,y$ is distinguished by any $u\in  B_j-M_j$.

Hence, $B$ is a local metric generator for $G[{\cal H}]$ and, as a consequence, 
  $$\dim_l(G[{\cal H}])\le \vert B \vert =\displaystyle\sum_{j\in J_{\cal H}}(\vert B_j\vert -\vert M_j\vert ) =
  \displaystyle\sum_{j\in J_{\cal H}}(\dim_l(G_j)-\alpha_j).
 $$
  Therefore, the result follows.
\end{proof}


\begin{theorem}\label{Minimalminimo}
Let $G[{\cal H}]$ be  a non-bipartite graph  obtained by point-attaching from a  family of connected graphs ${\cal H}=\{G_1, . . . ,G_k\}$. If for each $j\in [k]$ it holds that any minimal local metric generator for $G_j$ is minimum,   then  
$$\dim_l(G[{\cal H}])=\sum_{j\in J_{\cal H}}(\dim_l(G_j)-\alpha_j).$$
\end{theorem}

\begin{proof}
Since $G[{\cal H}]$ is a non-bipartite graph, any vertex belonging to a local metric basis of $G[{\cal H}]$ distinguishes  every pair of adjacent vertices included in a bipartite primary subgraph of $G[{\cal H}]$. Hence, we take a local metric basis $A$ of $G[{\cal H}]$  which does not contain vertices belonging to the bipartite primary subgraphs of $G[{\cal H}]$. \textit{i.e.}, for any $i\in [k]-J_{\cal H}$ it holds $A\cap V(G_i)=\emptyset$. Now, for each $j\in J_{\cal H}$ we define  $A_j=A\cap V(G_j)$.

  We claim  that $C_j \cup A_j$ is a local metric generator for $G_j$. Suppose that there exist two adjacent vertices $x,y\in V(G_j)$
 which are not distinguished by the elements of $ A_j$. In such a case,  there exists $x_r\in A_r$, $r\in J_{\cal H}-\{j\}$, which distinguishes  $x,y$, and so there must exists  $v\in C_j$ such that $G_r$ is a subgraph of $G_j(v^+)$ and, as a result, $v$ distinguishes the pair $x,y$. Hence,  $C_j \cup A_j$ is a local metric generator for $G_j$.
 
  Moreover, if $j\in J_{\cal H}$, then  for any   attachment vertex $w\in C_j$ it holds that $\vert A\cap V(G_j(w^+))\vert >0$, as $G_j(w^+)$ is not bipartite. Hence, given two adjacent vertices  $x,y\in V(G_j)$, which are distinguished by $w$, there exists $w'\in A_r\cap V(G_j(w^+)$, $r\in J_{\cal H}-\{j\}$, which distinguishes $x,y$, and so the minimality of $A$ leads to $C_j\cap A_j=\emptyset$. 
 
Now, if any minimal local metric generator for $G_j$ is minimum,   then 
  there exists  a set $C'_j\subseteq C_j$ such that $C'_j\cup A_j$ is a local metric basis for $G_j$.  
 Thus, $\vert C'_j \vert +\vert  A_j \vert=  \vert C'_j \cup A_j \vert = \dim_l(G_j)$. 
 Therefore,
 $$\dim_l(G[{\cal H}])=\vert A\vert =\sum_{j\in J_{\cal H}} \vert A_j\vert = \sum_{j\in J_{\cal H}}(\dim_l(G_j)-\vert C'_j \vert)\ge \sum_{j\in J_{\cal H}}(\dim_l(G_j)-\alpha_j).$$
 We conclude the proof by Theorem \ref{SecondTheoremSeparableGraphs}.
 \end{proof}


For any $j\in J_{{\cal H}}$ we define $\Gamma(G_j)$ as the family of local metric generators for $G_j$, and 
 $$\rho_j=\min_{S\subseteq V(G_j)}\left\{\vert S\vert :\; S\cup C_j\in \Gamma(G_j)\right\}.$$
 Also, any set for which the above minimum is attained will be denoted by $R_j$. Notice that  such a set  is not necessarily unique.

 With the above notation  in mind we can state our next result.

\begin{theorem}\label{mainTheoremSeparableGraphs}
For any   non-bipartite graph $G[{\cal H}]$ obtained by point-attaching from a family of connected graphs ${\cal H}=\{G_1, . . . ,G_k\}$,
$$\dim_l(G[{\cal H}])= \sum_{j\in J_{\cal H}}\rho_j.$$
\end{theorem}

\begin{proof}
We will show that  $X=\displaystyle\bigcup_{j\in J_{{\cal H}}}R_j$ 
is a local metric generator for $G[{\cal H}]$. 

First of all, note that by the structure of $G[{\cal H}]$ we have that for any $v\in C_j$, $j\in J_{{\cal H}}$,   there exists a non-bipartite primary subgraph
 $G_i$, which is a subgraph of  $G_j(v^+)$, such that $R_i\ne \emptyset$. To see this we 
  take a non-bipartite primary subgraph  $G_{j_1}$, which is a subgraph of $G_j(v^+)$, next, if $R_{j_1}=\emptyset$, then we take $v_1\in V(G_{j_1})-\{v\}$ and, as above, we take a non-bipartite primary subgraph  $G_{j_2}$, which is a subgraph of $G_j(v_{1}^+)$, and if $R_{j_2}=\emptyset$ then we repeat this process until obtain  a non-bipartite primary subgraph  $G_{j_t}$, which is a subgraph of 
    $G_j(v_{t-1}^+)$ such that  $R_{j_t}\ne \emptyset$ (at worst, we will arrive to a subgraph $G_j(v_{t-1}^+)$ containing only one non-bipartite primary subgraph). Hence, $X\ne \emptyset$ and, as a result, if $G_i$ is bipartite, then any pair of adjacent vertices $x,y\in V(G_i)$   is distinguished by any vertex belonging to $X$. 
    
    Now, if $x,y$ are adjacent in a non-bipartite primary subgraph $G_j$, then there exists $v\in R_j\cup C_j$ which distinguishes them. In the case that $v\in C_j$, we know that  there exists  a primary subgraph of  $G_j(v^+)$, such that $R_i\ne \emptyset$ and any $w\in R_i$ also distinguishes $x,y$. As a result, $X$ is a local metric generator for $G[{\cal H}]$. Therefore,
    $$\dim_l(G[{\cal H}])\le \vert X\vert = \sum_{j\in J_{\cal H}}\rho_j.$$

It remains to show that $\dim_l(G[{\cal H}])\ge \vert X\vert = \sum_{j\in J_{\cal H}}\rho_j.$
Since $G[{\cal H}]$ is a non-bipartite graph, any vertex belonging to a local metric basis of $G[{\cal H}]$ distinguishes  every pair of adjacent vertices included in a bipartite primary subgraph of $G[{\cal H}]$. Hence, we take a local metric basis $A$ of $G[{\cal H}]$  which does not contain vertices belonging to the bipartite primary subgraphs of $G[{\cal H}]$  \textit{ i.e}., for any $i\in [k]-J_{\cal H}$ it holds $A\cap V(G_i)=\emptyset$. For each $j\in J_{\cal H}$ we define  $A_j=A\cap V(G_j)$. Note that $A_j\cup C_j$ is a 
local metric generator for $G_j$ and, by the minimality of $A$, we have $A_j\cap C_j=\emptyset$. Hence,   $\vert A_j\vert\ge \vert  R_j \vert =\rho_j$. Therefore, 
$$\dim_l(G[{\cal H}])=\vert A\vert = \sum_{j\in J_{\cal H}} \vert A_j\vert\ge  \sum_{j\in J_{\cal H}}\rho_j.$$
 \end{proof}
 
If  $G_j$ is the only non-bipartite primary subgraph of $G[{\cal H}]$, then $\vert J_{\cal H}\vert=1$ and $\rho_j=\dim_l(G_j)$. Then we obtain the following particular case of Theorem 
\ref{mainTheoremSeparableGraphs}.

\begin{corollary}\label{OnlyOneNonbipartite}
Let $G[{\cal H}]$ be a graph  obtained by point-attaching from the family of connected graphs ${\cal H}=\{G_1, . . . ,G_k\}$. If  $G_j$ is the only non-bipartite primary subgraph of $G[{\cal H}]$, then 
$$\dim_l(G[{\cal H}])= \dim_l(G_j).$$
\end{corollary} 
 It is well-known that that a unicyclic graph $G$ is bipartite if and only if its cycle has even length. For the case of non-bipartite unicyclic graphs we can apply Corollary \ref{OnlyOneNonbipartite} to deduce that for any non-bipartite unicyclic graph $G$ it holds that $\dim_l(G)=2$. 
 
 There are other cases in which $\rho_j$ and $\alpha_j$ are very easy to obtain. For instance, if  $C_j=\{v\}$, then 
 $\rho_j=\dim_l(G_j)-\alpha_j$, where $\alpha_j=1$ if $v$ belongs to a local metric basis for $G_i$ and $\alpha_j=0$ in  
otherwise. Also, if $C_j=V(G_j)$, then $\rho_j=0$ and $\alpha_j=\dim_l(G_j)$.

The remain sections of this article are devoted to derive some consequences of Theorem \ref{mainTheoremSeparableGraphs}. We also give several families of graphs where the   equality of Theorem \ref{SecondTheoremSeparableGraphs} is achieved. 

\section{Rooted product graphs}
Rooted product graphs can be constructed as follows.  Let $\mathcal{H}$ be a
sequence of $n$ graphs $H_{1},H_{2},\ldots ,$ $H_{n}.$ In each of these graphs
a particular vertex $v_{i}$ is selected. This vertex will be called the root
of the graph $H_{i}.$ The \textit{rooted product graph} $G\circ \mathcal{H},$ is the
graph obtained by identifying the root of the graph $H_i$ with the $i$-th vertex
of $G,$ as defined by Godsil and Mckay \cite{Godsil1978}. Clearly, any rooted product graph is  is obtained by point-attaching from $G, H_1,H_2, . . . ,H_n$. Therefore,  as a consequence of Theorem \ref{mainTheoremSeparableGraphs} we obtain a formula for the local metric dimension of any rooted product graph. To begin with, we consider the case where every $H_i$ is a bipartite graph.

\begin{corollary}
Let $G$ be a connected graph of order $n\ge 2$ and let $\mathcal{H}$ be a
sequence of $n$ connected bipartite graphs $H_{1},H_{2},\ldots ,$ $H_{n}$. Then for any rooted product graph $G\circ \mathcal{H},$
$$ \dim_l(G\circ \mathcal{H})=\dim_l(G).$$
\end{corollary}

If every $H_i$ is non-bipartite, the result can be expressed as follows.

\begin{corollary}\label{mainResultRooted}
Let $G$ be a connected  graph of order $n\ge 2$ and let $\mathcal{H}$ be a
sequence of $n$ connected non-bipartite graphs $H_{1},H_{2},\ldots ,$ $H_{n}$. Then for any rooted product graph $G\circ \mathcal{H},$
$$ \dim_l(G\circ \mathcal{H})=\sum_{j=1}^n(\dim_l(H_j)-\alpha_j).$$
\end{corollary}

Note that in this case $\alpha_j=1$ if the root of $H_j$ belongs to a local metric basis of $H_j$ and $\alpha_j=0$ in otherwise.


Now we will restrict
ourselves to a particular case of rooted product graphs where the sequence
$H_{1},H_{2},\ldots ,$ $H_{n}$ consists of $n$ isomorphic graphs   of order $n'$, and will be
using in each of them the same root vertex $v.$ 
 The resulting 
rooted product graph is denoted by the expression $G\circ _{v}H$. In this case Corollary \ref{mainResultRooted} is simplified as follows.

\begin{remark}\label{Tnonb}
Let $H$ be a connected non-bipartite  graph and let $v$ be a vertex of $H.$

\begin{enumerate}[{\rm (i)}]
\item If $v$ does not belong to any metric basis for $H,$ then for any connected
graph $G$ of order $n,$
\[
\dim _{l}\left( G\circ _{v}H\right) =n\cdot \dim _{l}(H)
\]
\item If $v$ belongs to a  metric basis for $H,$ then for any connected graph $%
G $ of order $n\geq 2,$
\[
\dim _{l}\left( G\circ _{v}H\right) =n\cdot \left( \dim _{l}(H)-1\right)
\]

\end{enumerate}
\end{remark}


\begin{lemma}\label{Lemmale(n-3)}
If $H$ is a connected graph of order $n'$ with clique number $\omega
(H)=n'-1,$ and $G$ is a connected graph of order $n\ge 2$, then for any $v\in V(H),$
\[
\dim _{l}(G\circ_v H)=n(n'-3).
\]
\end{lemma}

\begin{proof}
Since $H$ has clique number $\omega (H)=n'-1$, by Theorem \ref{ldimb02} we have  $\dim_l(H)=n'-2$. To conclude the proof by Remark \ref{Tnonb} we need to prove that any vertex of  $H$ belongs to a local metric basis.
With this aim, we consider three vertices $v_i,v_j,v_k\in V(H)$ and a maximum clique $Q$ of $H$ such that  $v_i\not\in V(Q)$,   $v_j\in N_H(v_i)$ and $v_k\not\in N_H(v_i)$ (Here $N_H(x)$ denotes the set of neighbours that $x$ has in $H$). Then we have the following:
\begin{itemize}
\item Since $v_i$ distinguishes the pair of adjacent vertices $v_j,v_k$, the set   $B_i=V(H)-\{v_j,v_k\}$ is a local metric basis of $H$.
\item Since $v_iv_k\not\in E(H)$, the set,  $B_j=V(H)-\{v_i,v_k\}$ is a local metric basis of $H$.
\item Since $v_k$ distinguishes the pair of adjacent vertices $v_i,v_j$, the set  $B_k=V(H)-\{v_i,v_j\}$ is a local metric basis of $H$.
\end{itemize}
Therefore, any vertex of  $H$  belongs to a local metric basis.
\end{proof}

The equality $\dim _{l}(G\circ_v H)=n(n'-3)$  is not exclusive for connected graphs of order $n'$  with clique
number $\omega (H)=n'-1.$ Consider for instance the graph $H=\langle  v \rangle+(K_r\cup K_s)$, $r\ge 2$ and $s\ge 2$, \textit{i.e.,} $H$ is the graph $K_r\cup K_s$ together with all the edges joining  an isolated vertex $v$ to every vertex of  $K_r\cup K_s$.
In this case the order of $H$ is $n'=r+s+1,$ while its local metric dimension is $ \dim _{l}(H)=n'-3$.  Note however, that the vertex $v$ can not be in any local metric
basis. Hence, in this particular case for any connected graph $G$ of order $n\ge 2,$
the local metric dimension of the rooted product graph $G\circ _{v}H$ is calculated from Remark \ref{Tnonb}, giving
\[
\dim _{l}\left( G\circ _{v}H\right) =n\cdot \dim _{l}(H)=n(n'-3).
\]

\begin{proposition}\label{theoremBounds}
Let $G$ be a connected graph of order $n\ge 2$.  Let $H$ be a connected non-bipartite graph of order $n'$ and let  $v\in V(H)$. Then the following assertions hold.

\begin{enumerate}[{\rm(i)}]
\item  $n\leq \dim _{l}(G\circ _{v}H)\leq n(n'-2).$

\item $\dim _{l}(G\circ _{v}H)=n$ if and only if $\dim_l(H)=2$ and the root vertex $v$ belongs to any local metric basis of $H$.

\item $\dim _{l}(G\circ _{v}H)=n(n'-2)$ if and only if  $H \cong K_{n'}$.

\item If  $H \not\cong K_{n'}$, then  $\dim _{l}(G\circ _{v}H)\leq n(n'-3).$
\end{enumerate}
\end{proposition}

\begin{proof}
Remark \ref{Tnonb}  directly leads to the lower bound. Note that $\dim_l(H)\ge 2$, as $H$ is not bipartite. Now, if $v$ belongs to a local metric basis of $H$ and $\dim_l(H)=2$, then Remark \ref{Tnonb} (ii) leads to $\dim _{l}(G\circ _{v}H)= n$. Otherwise, 
 if $v$ does not belong to any local metric basis of $H$, then Remark \ref{Tnonb} leads to $\dim _{l}(G\circ _{v}H)\ge 2n$. This proves (ii).

Now, if $H\cong  K_{n'}$, then $\dim_l(H)= n'-1$ and, since $v$ belongs to a local metric basis of $H$, Remark \ref{Tnonb} (ii) leads to $\dim _{l}(G\circ _{v}H)= n(n'-2)$. On the other hand, if $H$ is a connected non-complete graph   of order $n'$, then we have $\dim_l(H)\le n'-2$. 
So, Remark \ref{Tnonb}  leads to the upper bound. 

Note that if  $\dim_l(H)=n'-2$, then  Theorem \ref{ldimb02} and Lemma \ref{Lemmale(n-3)} lead to $\dim_l(G\circ _{v}H)\le n(n'-3)$. 
Thus, (iii) and (iv) follows. 
\end{proof}

\section{Corona product graphs}

Let $G$  be a graphs of order $n$ and let ${\cal H}=\{H_1,H_2,...,H_n\}$ be a family of graphs. Recall that the corona product $G\odot {\cal H}$ is defined as the graph obtained from $G$ and ${\cal H}$ by taking one copy of $G$ and joining by an edge each vertex from $H_i$ with the $i^{th}$-vertex of $G$, \cite{Frucht1970}. The join $G+H$ is defined as the graph obtained from disjoint graphs $G$ and $H$ by taking one copy of $G$ and one copy of $H$ and joining by an edge each vertex of $G$ with each vertex of $H$. Notice that the particular case of corona graph $K_1\odot H$ is isomorphic to the join graph $K_1+H$. 
We can obtain any corona graph $G\odot {\cal H}$ by point-attaching from $G, K_1+H_1,K_1+H_2, . . . ,K_1+H_n$. Note that if $H_i$ is a non-trivial graph, then the primary subgraph $K_1+H_i$ is not bipartite. In fact, we can see the corona graph as a particular case of rooted product graph.

\begin{corollary}\label{generalResultCorona}
Let $G$ be a connected  graph of order $n\ge 2$ and let $\mathcal{H}$ be a
sequence of $n$ non-empty graphs $H_{1},H_{2},\ldots ,$ $H_{n}$. Then for any corona product graph $G\circ \mathcal{H},$
$$ \dim_l(G\odot \mathcal{H})=\sum_{j=1}^n(\dim_l(K_1+H_j)-\alpha_j).$$
\end{corollary}

Note that in this case $\alpha_j=1$ if the vertex of $K_1$ belongs to a local metric basis of $K_1+H_j$ and $\alpha_j=0$ in otherwise.


The particular case of corona product graphs where the sequence
$H_{1},H_{2},\ldots ,$ $H_{n}$ consists of $n$ isomorphic graphs  of order $n'$ was previously studied in \cite{Rodriguez-Velazquez2013LDimCorona,Rodriguez-Velazquez-Fernau2013}.  
 The resulting  corona graph is denoted by the expression $G\odot H$. 
As a  particular case of Corollary \ref{generalResultCorona} we derive the next result which was previously obtained in \cite{Rodriguez-Velazquez2013LDimCorona}.

\begin{remark}{\rm \cite{Rodriguez-Velazquez2013LDimCorona}} \label{mainTheorem}
Let $H$ be a non-empty graph. The following assertions hold.
\begin{enumerate}[{\rm (i)}]
\item  If the vertex of $K_1$ does not belong to any local metric basis for $K_1+H$, then for any connected graph $G$ of order $n$,
$$\dim_l(G\odot H)=n\cdot \dim_l(K_1+H).$$

\item   If the vertex of $K_1$  belongs to a local metric basis for $K_1+H$, then for any connected graph $G$ of order $n\ge 2$,
$$\dim_l(G\odot H)=n  (\dim_l(K_1+H)-1).$$
\end{enumerate}
\end{remark}

The reader is referred to \cite{Rodriguez-Velazquez2013LDimCorona,Rodriguez-Velazquez-Fernau2013}  for a moire detailed study on the local metric dimension of corona product graphs.

\section{Block graphs}
A \textit{block graph} is a graph whose blocks are cliques. Since any block graph is obtained by point-attaching from $G_1=K_{t_1},G_2=K_{t_2}, . . . ,G_k=K_{t_k}$,  as a consequence of Theorem \ref{mainTheoremSeparableGraphs} we obtain a formula for the local metric dimension of any block graph. Our next result shows how the formula is reduced when every block has order $t_i\ge 3$.

\begin{corollary}
Let ${\cal H}=\{G_1=K_{t_1},G_2=K_{t_2}, . . . ,G_k=K_{t_k}\}$ be a finite sequence of pairwise disjoint complete  graphs of order $t_i\ge 3$, $i=1,...,k.$ Then for any block graph $G[{\cal H}]$,
$$\dim_l(G[{\cal H}])  =\sum_{j=1}^k(t_j-1-\alpha_j).$$
\end{corollary}
   
In this case $\alpha_j$ becomes $t_j-1$  if every vertex of $K_{t_j}$ is a cut vertex of $G[{\cal H}]$  and it becomes the number of cut vertices of $G[{\cal H}]$ belonging to the clique $K_{t_i}$ in otherwise.

\section{Bouquet of graphs}
Let ${\cal H}=\{G_1,G_2, . . . ,G_k\}$ be a finite sequence of pairwise disjoint connected graphs and let
$x_i \in V(G_i)$. By definition, the \textit{bouquet} ${\cal H}_x$ of the graphs in ${\cal H}$ with respect to the
vertices $\{x_i\}_{i=1}^k$ is obtained by identifying the vertices $x_1, x_2, . . . , x_k$ with a new vertex $x$.
Clearly, the bouquet ${\cal H}_x$ is a graph obtained by point-attaching from $G_1,G_2, . . . ,G_k$. Therefore,  as a consequence of Theorem \ref{mainTheoremSeparableGraphs} we obtain the following result. 

\begin{corollary}
Let ${\cal H}=\{G_1,G_2, . . . ,G_k\}$ be a finite sequence of pairwise disjoint connected graphs and let
$x_i \in V(G_i)$ such that $J_{\cal H}\ne \emptyset$. If ${\cal H}_x$ is the bouquet obtained from ${\cal H}$ by identifying the vertices $x_1, x_2, . . . , x_k$ with a new vertex $x$, then
$$\dim_l({\cal H}_x)  =\sum_{j\in J_{\cal H}}(\dim_l(G_j)-\alpha_j).$$
\end{corollary}

Note that in this case $\alpha_i=1$ if $x_i$ belongs to a local metric basis of $G_i$ and $\alpha_i=0$ in otherwise.

\section{Chain of graphs}
Let ${\cal H}=\{G_1,G_2, . . . ,G_k\}$ be a finite sequence of pairwise disjoint connected non-trivial graphs and let $x_i,y_i\in V(G_i)$.
 By definition, the \textit{chain} ${\cal C( H)}$ of the graphs in ${\cal H}$ with respect to the set of vertices $\{y_1,x_k\}\cup \left( \cup _{i=2}^{k-1} \{x_i,y_i\} \right)$ is the connected graph obtained by identifying the vertex $y_i$ with the vertex $x_{i+1}$ for $i\in [k-1]$. Clearly, the chain ${\cal C(H)}$ is a graph obtained by point-attaching from $G_1,G_2, . . . ,G_k$. 

\begin{figure}[h]
   \centering
\begin{tikzpicture}
\draw  (-2,0) ellipse (2cm and 0.8cm);
\filldraw[fill opacity=1,fill=black]  (0,0) circle (0.12cm);
\node at (0,-0.8) {$y_1=x_2$ };
\draw  (2,0) ellipse (2cm and 0.8cm);
\draw  (6,0) ellipse (2cm and 0.8cm);
\filldraw[fill opacity=1,fill=black]  (4,0) circle (0.12cm);
\node at (4,-0.8) {$y_2=x_3$ };
\draw  (10,0) ellipse (2cm and 0.8cm);
\filldraw[fill opacity=1,fill=black]  (8,0) circle (0.12cm);
\node at (8,-0.8) {$y_3=x_4$ };
\node at (-2,0) {$G_1$ };
\node at (2,0) {$G_2$ };
\node at (6,0) {$G_3$ };
\node at (10,0) {$G_4$ };
\end{tikzpicture}
\caption{A chain ${\cal C( H)}$ obtained by point-attaching from ${\cal H}=\{G_1,G_2, G_3,G_4\}$.}
\label{point-attaching-Chain}
\end{figure}
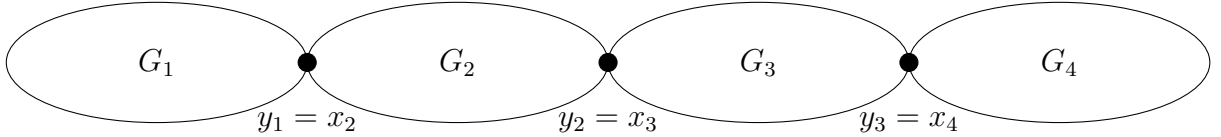 
 
  For every $j\in J_{\cal H}$ we say that $x_j$ is \textit{replaceable} in ${\cal C(H)}$ if and only if there exists a local metric basis $B_j$ of $G_j$ such that $x_j\in B_j$
and there exists $k < j$ such that $G_k$ is a non-bipartite primary graph.
Analogously, we say that $y_j$ is \textit{replaceable} in ${\cal C(H)}$  if and only if there exists a local metric basis $B'_j$ of $G_j$ such that $y_j\in B'_j$
and there exists $k > j$ such that $G_k$ is a non-bipartite primary subgraph. We say that $x_j$ and $y_j$ are \textit{simultaneously replaceable} in ${\cal C(H)}$ if both are replaceable in ${\cal C(H)}$ and there exists  a local metric basis  of $G_j$ containing both $x_j$ and $y_j$.

The formula for $\dim_l({\cal C(H)})$ is directly obtained from Theorem \ref{mainTheoremSeparableGraphs}. In this case we have the following  possibilities for the value of $\rho_j$. 

\begin{itemize}
 \item If $1\in J_{\cal H}$ and  $y_1$ is replaceable in $\cal C(H)$, then $\rho_1=\dim_l(G_1)-1$.
 
 \item If $1\in J_{\cal H}$ and $y_1$ is not replaceable in $\cal C(H)$, then $\rho_1=\dim_l(G_1)$.

 \item If  $k\in J_{\cal H}$ and $x_k$ is replaceable in $\cal C(H)$, then $\rho_k=\dim_l(G_1)-1$.
 
 \item If $k\in J_{\cal H}$ and  $x_k$ is not replaceable in $\cal C(H)$, then $\rho_k=\dim_l(G_1)$.   
  \end{itemize}
For $j\in  J_{\cal H} \cap \{2,..., k-1\}$ we have the following possibilities.
\begin{itemize}
\item If neither $x_j$ nor $y_j$ is replaceable in $\cal C(H)$, then either $\rho_j=\dim_l(G_j)$ or $\rho_j=\dim_l(G_j)-1$.

\item If $x_j$ and $y_j$ are simultaneously replaceable in $\cal C(H)$, then $\rho_j=\dim_l(G_j)-2$.

\item If $x_j$ and $y_j$ are not simultaneously replaceable in $\cal C(H)$ and $x_j$ (or $y_j$) is replaceable in $\cal C(H)$, then $\rho_j=\dim_l(G_j)-1$.
\end{itemize}

\end{document}